\documentclass[11pt]{Mytran-l}
\usepackage[ansinew]{inputenc}
\usepackage[centertags]{amsmath}
\usepackage{exscale}
\usepackage{relsize}
\usepackage{amsfonts}
\usepackage{amssymb}
\usepackage{amsthm}
\usepackage{newlfont}
\usepackage{color}
\usepackage[colorlinks,citecolor=blue,urlcolor=blue,breaklinks=true]{hyperref}

\mathchardef\ordinarycolon\mathcode`\:
  \mathcode`\:=\string"8000
  \begingroup \catcode`\:=\active
    \gdef:{\mathrel{\mathop\ordinarycolon}}
  \endgroup
\def\N{{\Bbb N}}

\newtheorem*{mthm*}{Main Theorem}

\newtheorem{thm}{Theorem}

\newtheorem{prop}[thm]{Proposition}

\newtheorem{rem}[thm]{Remark}
\newtheorem{Cor}[thm]{Corollary}

\newtheorem{conj}{Conjecture}

\newcommand{\D}{\mathcal{D}}
\newcommand{\B}{\mathcal{B}}
\newcommand{\aut}{\emph{Aut}}
\newcommand{\Aut}{\mbox{Aut}}

\newcommand{\G}{\mathit{\Gamma}}
\newcommand{\Si}{\mathit{\Sigma}}

\begin{document}
\title[The Cameron-Praeger Conjecture]
{On the Cameron-Praeger Conjecture}

\author{Michael Huber}



\date{November 21, 2008; and in revised form April 18, 2009}

\commby{The Managing Editor}

\maketitle
\vspace*{-0.5cm}

\begin{center}
{\footnotesize Wilhelm-Schickard-Institute for Computer Science\\
University of Tuebingen\\
Sand~13, D-72076 Tuebingen, Germany\\
E-mail: {\texttt{michael.huber@uni-tuebingen.de}}}
\end{center}

\smallskip

\begin{abstract}
This paper takes a significant step towards confirming a long-standing and far-reaching conjecture of Peter~J.~Cameron and
\linebreak Cheryl~E.~Praeger. They conjectured in 1993 that there are no non-trivial block-transitive \mbox{$6$-designs}. We prove that
the Cameron-Praeger conjecture is true for the important case of non-trivial Steiner \mbox{$6$-designs}, i.e.
for \mbox{$6$-$(v,k,\lambda)$} designs with \mbox{$\lambda=1$}, except possibly when
the group is \mbox{$P \mathit{\Gamma} L(2,p^e)$} with $p=2$ or $3$, and $e$ is an odd prime power.
\end{abstract}

\smallskip


\section{Introduction}\label{intro}

The characterization of combinatorial or geometric
structures in terms of their groups of automorphisms has attracted
considerable interest in the last decades and is now commonly viewed
as a natural generalization of Felix Klein’s Erlangen program
(1872).
There has been recent progress in particular on the characterization
of Steiner \mbox{$t$-designs} which admit groups of automorphisms
with sufficiently strong symmetry properties: The author classified
all flag-transitive Steiner \mbox{$t$-designs} with $t>2$
(see~\cite{Hu2001,Hu2005,Hu_Habil2005,Hu2007,Hu2007a}
and~\cite{Hu2008} for a monograph). In particular, he showed
in~\cite{Hu2007} that no non-trivial flag-transitive Steiner
\mbox{$6$-design} can exist. These results answer a series of
40-year-old problems and generalize theorems of
J.~Tits~\cite{Tits1964} and H.~Lüneburg~\cite{Luene1965}. Previously,
F.~Buekenhout, A.~Delandtsheer, J.~Doyen, P.~Kleidman, M.~Liebeck,
and J.~Saxl~\cite{Buek1990,Del2001,Kleid1990,Lieb1998,Saxl2002} had
characterized all flag-transitive Steiner
\mbox{2-designs}, up to the \mbox{$1$-dimensional} affine case. All these classification results rely on the
classification of the finite simple groups.

In 1993, P.~J.~Cameron and C.~E.~Praeger (\cite[Conj.\,1.2]{CamPrae1993}) conjectured that there
are no non-trivial block-transitive \mbox{$6$-designs}.\footnote{see
also Kourovka Notebook~\cite[Problem\,11.45]{Kor1992}, and Peter Cameron's conjectures online at
\url{http://www.maths.qmw.ac.uk/~pjc/cameronconjs.html}.}
Our main result is as follows:

\begin{mthm*}\label{mainthm}
Let $\mathcal{D}=(X,\mathcal{B},I)$ be a non-trivial Steiner
\mbox{$6$-design}. Then \mbox{$G \leq \emph{Aut}(\mathcal{D})$}
cannot act block-transitively on $\mathcal{D}$, except possibly when
\linebreak \mbox{$G= P \mathit{\Gamma} L(2,p^e)$} with $p=2$ or $3$
and $e$ is an odd prime power.
\end{mthm*}

The result has been announced (without proof) in a recent
paper~\cite{Hu_mmics2008} on the existence problem for Steiner
\mbox{$t$-designs} for large values of $t$. The proof makes use of
the classification of the finite \mbox{$3$-homogeneous} permutation
groups, which in turn relies on the classification of the finite
simple groups. It will be given in Section~\ref{proof}. Preliminary
results which are important for the remainder of the paper are
collected in Section~\ref{Prelim}.

\medskip


\section{Definitions and Notations}\label{Defn}

For positive integers $t \leq k \leq v$ and $\lambda$,
we define a \mbox{\emph{$t$-$(v,k,\lambda)$ design}}
to be a finite incidence structure
\mbox{$\mathcal{D}=(X,\mathcal{B},I)$}, where $X$ denotes a set of
\emph{points}, $\left| X \right| =v$, and $\mathcal{B}$ a set of
\emph{blocks}, $\left| \mathcal{B} \right| =b$, with the following
regularity properties: each block $B \in \mathcal{B}$ is incident
with $k$ points, and each \mbox{$t$-subset} of $X$ is incident with
$\lambda$ blocks. A \emph{flag} of $\mathcal{D}$ is an incident
point-block pair $(x,B) \in I$ with $x \in X$ and $B \in
\mathcal{B}$.

For historical reasons, a \mbox{$t$-$(v,k,\lambda)$ design} with
$\lambda =1$ is called a \emph{Steiner \mbox{$t$-design}} (sometimes
also a \emph{Steiner system}). We note that in this case each block
is determined by the set of points which are incident with it, and
thus can be identified with a \mbox{$k$-subset} of $X$ in a unique
way. If $t<k<v$, then we speak of a \emph{non-trivial} Steiner
\mbox{$t$-design}. There are many infinite classes of
Steiner \mbox{$t$-designs} for $t=2$ and $3$, however for $t=4$ and
$5$ only a finite number are known. For a detailed treatment of
combinatorial designs, we refer
to~\cite{BJL1999,crc06,hall86,hupi85,stin04}. In
particular,~\cite{BJL1999,crc06} provide encyclopedic accounts of
key results and contain existence tables with known parameter sets.

In what follows, we are interested in \mbox{$t$-designs} which admit
groups of automorphisms with sufficiently strong symmetry properties
such as transitivity on the blocks or on the flags. We consider
automorphisms of a \mbox{$t$-design} $\mathcal{D}$ as pairs of
permutations on $X$ and $\mathcal{B}$ which preserve incidence, and
call a group \mbox{$G \leq \mbox{Aut} (\mathcal{D})$} of
automorphisms of $\mathcal{D}$ \emph{block-transitive} (respectively
\emph{flag-transitive}, \emph{point \mbox{$t$-transitive}},
\emph{point \mbox{$t$-homogeneous}}) if $G$ acts transitively on the
blocks (respectively transitively on the flags,
\mbox{$t$-transitively} on the points, \mbox{$t$-homogeneously} on
the points) of $\mathcal{D}$. For short, $\mathcal{D}$ is said to
be, e.g., block-transitive if $\mathcal{D}$ admits a
block-transitive group of automorphisms.

For \mbox{$\mathcal{D}=(X,\mathcal{B},I)$} a Steiner
\mbox{$t$-design} with \mbox{$G \leq \mbox{Aut} (\mathcal{D})$}, let
$G_x$ denote the stabilizer of a point $x \in X$, and $G_B$ the
setwise stabilizer of a block $B \in \mathcal{B}$. For $x, y \in X$
and $B \in \mathcal{B}$, we define $G_{xy}= G_x \cap G_y$.

\medskip


\section{Preliminary Results}\label{Prelim}

\subsection{Combinatorial Results}\hfill

\smallskip

Basic necessary conditions for the existence of \mbox{$t$-designs}
can be obtained via elementary counting arguments (see, for
instance,~\cite{BJL1999}):

\begin{prop}\label{s-design}
Let $\mathcal{D}=(X,\mathcal{B},I)$ be a \mbox{$t$-$(v,k,\lambda)$}
design, and for a positive integer $s \leq t$, let $S \subseteq X$
with $\left|S\right|=s$. Then the total number of blocks incident
with each element of $S$ is given by
\[\lambda_s = \lambda \frac{{v-s \choose t-s}}{{k-s \choose t-s}}.\]
In particular, for $t\geq 2$, a \mbox{$t$-$(v,k,\lambda)$} design is
also an \mbox{$s$-$(v,k,\lambda_s)$} design.
\end{prop}

It is customary to set $r:= \lambda_1$ denoting the total
number of blocks incident with a given point.

\begin{Cor}\label{Comb_t=5}
Let $\mathcal{D}=(X,\mathcal{B},I)$ be a \mbox{$t$-$(v,k,\lambda)$}
design. Then the following holds:
\begin{enumerate}

\item[{(a)}] $bk = vr.$

\smallskip

\item[{(b)}] $\displaystyle{{v \choose t} \lambda = b {k \choose t}.}$

\smallskip

\item[{(c)}] $r(k-1)=\lambda_2(v-1)$ for $t \geq 2$.

\end{enumerate}
\end{Cor}

\begin{Cor}\label{divCond}
Let $\mathcal{D}=(X,\mathcal{B},I)$ be a \mbox{$t$-$(v,k,\lambda)$}
design. Then
\[\lambda {v-s \choose t-s} \equiv \, 0\; \emph{(mod}\;\, {k-s \choose t-s})\]
for each positive integer $s \leq t$.
\end{Cor}

For non-trivial Steiner \mbox{$t$-designs} lower bounds for $v$ in
terms of $k$ and $t$ can be given (see
P.~Cameron~\cite[Thm.\,3A.4]{Cam1976}, and
J.~Tits~\cite[Prop.\,2.2]{Tits1964}):

\begin{prop}\label{Cam}
If $\mathcal{D}=(X,\mathcal{B},I)$ is a non-trivial Steiner
\mbox{$t$-design}, then the following holds:
\begin{enumerate}

\item[{(a)}] \emph{(Tits~1964):} \hspace{0.2cm} $v\geq (t+1)(k-t+1).$

\smallskip

\item[{(b)}] \emph{(Cameron~1976):} \hspace{0.2cm} \mbox{$v-t+1 \geq (k-t+2)(k-t+1)$} for $t>2$. If equality
holds, then
\smallskip
$(t,k,v)=(3,4,8),(3,6,22),(3,12,112),(4,7,23)$, or $(5,8,24)$.
\end{enumerate}
\end{prop}

In the case when \mbox{$t=6$}, we deduce from Part~(b) the following upper bound for the positive integer
$k$.

\begin{Cor}\label{Cameron_t=5}
Let $\D=(X,\B,I)$ be a non-trivial Steiner \mbox{$t$-design} with
\mbox{$t=6$}. Then
\[k \leq \Bigl\lfloor \sqrt{v-\textstyle{\frac{19}{4}}} + \textstyle{\frac{9}{2}} \Bigr\rfloor.\]
\end{Cor}

\smallskip


\subsection{Highly Symmetric Designs}\hfill

\smallskip

We will now focus on \mbox{$t$-designs} which admit groups
of automorphisms with sufficiently strong symmetry properties. One
of the reasons for this consideration of highly symmetric designs is
a general view that, while the existence of combinatorial objects is
of interest, they are even more fascinating when they have a rich
group of symmetries.

One of the early important results regarding highly symmetric
designs is due to R.~Block~\cite[Thm.\,2]{Block1965}:

\begin{prop}{\em (Block~1965).}\label{BlocksLemma2}
Let $\mathcal{D}=(X,\mathcal{B},I)$ be a non-trivial
\mbox{$t$-$(v,k,\lambda)$} design with $t \geq 2$. If $G \leq
\emph{Aut}(\mathcal{D})$ acts block-transitively on $\mathcal{D}$,
then $G$ acts point-transitively on $\mathcal{D}$.
\end{prop}

For a \mbox{$2$-$(v,k,1)$} design $\mathcal{D}$, it is elementary
that the point \mbox{$2$-transitivity} of \mbox{$G \leq
\mbox{Aut}(\mathcal{D})$} implies its flag-transitivity. For
\mbox{$2$-$(v,k,\lambda)$} designs, this implication remains true if
$r$ and $\lambda$ are relatively prime
(cf.~\cite[Chap.\,2.3,\,Lemma\,8]{demb68}). However, for
\mbox{$t$-$(v,k,\lambda)$} designs with $t \geq 3$, it can be
deduced from Proposition~\ref{BlocksLemma2} that always the converse
holds (see~\cite{Buek1968} or~\cite[Lemma\,2]{Hu2001}):

\begin{prop}\label{flag2trs}
Let $\mathcal{D}=(X,\mathcal{B},I)$ be a non-trivial
\mbox{$t$-$(v,k,\lambda)$} design with \mbox{$t \geq 3$}. If
\mbox{$G \leq \emph{Aut}(\mathcal{D})$} acts flag-transitively on
$\mathcal{D}$, then $G$ acts point \linebreak
\mbox{$2$-transitively} on $\mathcal{D}$.
\end{prop}

Investigating highly symmetric \mbox{$t$-designs} for large values
of $t$, P.~Cameron and C.~Praeger~\cite[Thm.\,2.1]{CamPrae1993}
derived from Proposition~\ref{BlocksLemma2} and a combinatorial result of
D.~Ray-Chaudhuri and R.~Wilson~\cite[Thm.\,1]{Ray-ChWil1975}
the following assertion:

\begin{prop}{\em (Cameron \& Praeger~1993).}\label{flag3hom}
Let $\mathcal{D}=(X,\mathcal{B},I)$ be a \mbox{$t$-$(v,k,\lambda)$}
design with $t\geq 2$. Then, the following holds:

\begin{enumerate}

\item[{(a)}] If \mbox{$G \leq \emph{Aut}(\mathcal{D})$} acts block-transitively on $\mathcal{D}$,
then $G$ also acts point \linebreak \mbox{$\lfloor t/2
\rfloor$-homogeneously} on $\mathcal{D}$.

\smallskip

\item[{(b)}] If \mbox{$G \leq \emph{Aut}(\mathcal{D})$} acts flag-transitively on $\mathcal{D}$,
then $G$ also acts point \linebreak \mbox{$\lfloor (t+1)/2
\rfloor$-homogeneously} on $\mathcal{D}$.

\end{enumerate}
\end{prop}

As for $t \geq 7$ the flag-transitivity, respectively for $t \geq 8$
the block-transitivity of \mbox{$G \leq \mbox{Aut} (\mathcal{D})$}
implies at least its point \mbox{$4$-homogeneity}, they obtained the
following restrictions as a consequence of the finite simple group
classification (cf.~\cite[Thm.\,1.1]{CamPrae1993}):

\begin{thm}{\em (Cameron \& Praeger~1993).}
Let $\mathcal{D}=(X,\mathcal{B},I)$ be a \mbox{$t$-$(v,k,\lambda)$}
design. If \mbox{$G \leq \emph{Aut}(\mathcal{D})$} acts
block-transitively on $\mathcal{D}$ then $t \leq 7$, while if
\mbox{$G \leq \emph{Aut}(\mathcal{D})$} acts flag-transitively on
$\mathcal{D}$ then $t \leq 6$.
\end{thm}

Moreover, they formulated the following far-reaching conjecture
(cf.~\cite[Conj.\,1.2]{CamPrae1993}):

\begin{conj}{\em (Cameron \& Praeger~1993).}
There are no non-trivial block-transitive \mbox{$6$-designs}.
\end{conj}

\smallskip


\subsection{Finite $3$-homogeneous Permutation Groups}\hfill

\smallskip

In order to investigate all block-transitive Steiner \mbox{$6$-designs}, we can as a consequence
of Proposition~\ref{flag3hom}~(a) make use of the classification of
all finite \mbox{$3$-homogeneous} permutation groups, which itself
relies on the classification of all finite simple groups
(cf.~\cite{Cam1981,Gor1982,Kant1972,Lieb1987,LivWag1965}).

Let $G$ be a finite \mbox{$3$-homogeneous} permutation group on a
set $X$ with \mbox{$\left|X\right| \geq 4$}. Then $G$ is either of

{\bf (A) Affine Type:} $G$ contains a regular normal subgroup $T$
which is elementary Abelian of order $v=2^d$. If we identify $G$
with a group of affine transformations
\[x \mapsto x^g+u\]
of $V=V(d,2)$, where $g \in G_0$ and $u \in V$, then
one of the following occurs:

\begin{enumerate}

\smallskip

\item[(1)] $G \cong AGL(1,8)$, $A \mathit{\Gamma} L(1,8)$, or $A \mathit{\Gamma} L(1,32)$

\smallskip

\item[(2)] $G_0 \cong SL(d,2)$, $d \geq 2$

\smallskip

\item[(3)] $G_0 \cong A_7$, $v=2^4$

\end{enumerate}

\smallskip

or

\medskip

{\bf (B) Almost Simple Type:} $G$ contains a simple normal subgroup
$N$, and \mbox{$N \leq G \leq \Aut(N)$}. In particular, one of the
following holds, where $N$ and $v=|X|$ are given as follows:
\begin{enumerate}

\smallskip

\item[(1)] $A_v$, $v \geq 5$

\smallskip

\item[(2)] $PSL(2,q)$, $q>3$, $v=q+1$

\smallskip

\item[(3)] $M_v$, $v=11,12,22,23,24$ \hfill (Mathieu groups)

\smallskip

\item[(4)] $M_{11}$, $v=12$

\end{enumerate}

\medskip

We note that if $q$ is odd, then $PSL(2,q)$ is $3$-homogeneous for
\mbox{$q \equiv 3$ (mod $4$)}, but not for \mbox{$q \equiv 1$ (mod
$4$)}, and hence not every group $G$ of almost simple type
satisfying (2) is $3$-homogeneous on $X$. For required basic
properties of the listed groups, we refer, e.g.,
to~\cite{Atlas1985},~\cite{HupI1967},~\cite[Ch.\,2,\,5]{KlLi1990}.

\smallskip

\begin{rem} \label{equa_t=5}
\emph{If \mbox{$G \leq \mbox{Aut}(\mathcal{D})$} acts
block-transitively on any Steiner \mbox{$t$-design} $\mathcal{D}$
with $t\geq 6$, then by Proposition~\ref{flag3hom}~(a), $G$ acts
point \mbox{$3$-homogeneously} and in particular point
\mbox{$2$-transitively} on $\mathcal{D}$. Applying
Corollary~\ref{Comb_t=5}~(b) yields the equation
\[b=\frac{{v \choose t}}{{k \choose
t}}=\frac{v(v-1) \left|G_{xy}\right|}{\left| G_B \right|},\] where
$x$ and $y$ are two distinct points in $X$ and $B$ is a block in
$\mathcal{B}$.}
\end{rem}

\medskip


\section{Proof of the Main Theorem}\label{proof}

\emph{Let $\D=(X,\B,I)$ be a non-trivial Steiner \mbox{$6$-design} with \mbox{$G \leq \aut(\D)$}
acting block-transitively on $\D$ throughout the proof}. We recall
that due to Proposition~\ref{flag3hom}~(a), we may restrict
ourselves to the consideration of the finite \mbox{$3$-homogeneous}
permutation groups listed in Section~\ref{Prelim}.
Clearly, in the following we may assume that $k>6$ as trivial
Steiner \mbox{$6$-designs} are excluded.

\smallskip


\subsection{Groups of Automorphisms of Affine Type}\hfill

\medskip

\emph{Case} (1): $G \cong AGL(1,8)$, $A \mathit{\Gamma} L(1,8)$, or $A \mathit{\Gamma} L(1,32)$.

\smallskip

If \mbox{$v=8$}, then Corollary~\ref{Cameron_t=5} yields
$k \leq 6$, a contradiction. For $v=32$, Corollary~\ref{Cameron_t=5} implies that $k=7$, $8$ or $9$; for each of these values, $29$ divides $b$, and
so divides $\left|G\right|$ by block-transitivity, a contradiction since $29$ does not divide $\left| A \mathit{\Gamma} L(1,32) \right|$.

\medskip

\emph{Case} (2): $G_0 \cong SL(d,2)$, $d \geq 2$.

\smallskip

Here $v=2^d > k > 6$. For $d=3$, we have $v=8$, already ruled out in Case~(1).
So, we may assume that $d>3$. Any six distinct points being non-coplanar in
$AG(d,2)$, they generate an affine subspace of dimension at
least $3$. Let $\mathcal{E}$ be the
\mbox{$3$-dimensional} vector subspace spanned by the first three basis vectors $e_1,e_2,e_3$
of the vector space $V=V(d,2)$.
Then the pointwise stabilizer of $\mathcal{E}$ in $SL(d,2)$ (and therefore also in $G$)
acts point-transitively on \mbox{$V \setminus \mathcal{E}$}. If the
unique block $B \in \B$ which is incident with the \mbox{$6$-subset}
\mbox{$\{0,e_1,e_2,e_3,e_1+e_2,e_2+e_3\}$} contains some point
outside $\mathcal{E}$, then $B$ contains all points of
\mbox{$V \setminus \mathcal{E}$}, and so \mbox{$k
\geq v-2$}, a contradiction to
Corollary~\ref{Cameron_t=5}. Hence $B$ lies completely in
$\mathcal{E}$, and so $k \leq 8$. On the other hand, for $\D$
to be a block-transitive \mbox{$6$-design} admitting \mbox{$G \leq
\Aut(\D)$}, we deduce from~\cite[Prop.\,3.6\,(b)]{CamPrae1993} the
necessary condition that $2^d-3$ must divide $k \choose 4$, and
hence it follows for each respective value of $k$ that $d=3$,
contradicting our assumption.

\medskip

\emph{Case} (3): $G_0 \cong A_7$, $v=2^4$.

\smallskip

For $v=2^4$, we have $k \leq 7$ by Corollary~\ref{Cameron_t=5}, contradicting Proposition~\ref{s-design} since $r=\lambda_1$ is not an integer.

\smallskip


\subsection{\mbox{Groups of Automorphisms of Almost Simple Type}}\hfill

\medskip
\emph{Case} (1): $N=A_v$, $v \geq 5$.
\medskip

Since $\mathcal{D}$ is non-trivial with $k>6$, we may assume that $v \geq 8$. Then $A_v$, hence also $G$,
is \mbox{$6$-transitive} on $X$, and so cannot act on any non-trivial Steiner
\mbox{$6$-design} by~\cite[Thm.\,3]{Kant1985}.

\medskip

\emph{Case} (2): $N=PSL(2,q)$, $v=q+1$, $q=p^e >3$.

\smallskip

Here $\Aut(N)= P \mathit{\Gamma} L (2,q)$, and $\left| G \right| =
(q+1)q \frac{(q-1)}{n}a$ with $n=(2,q-1)$ and $a \mid ne$. We may again
assume that $v=q+1 \geq 8$.

\smallskip

\emph{We will first assume that $N=G$.} Then, by
Remark~\ref{equa_t=5}, we obtain
\begin{equation}\label{Eq-0}
(q-2)(q-3)(q-4) \left| PSL (2,q)_{B} \right| n  =
k(k-1)(k-2)(k-3)(k-4)(k-5).
\end{equation}
In view of Proposition~\ref{Cam}~(b), we have
\begin{equation}\label{Eq-A}
q-4 \geq (k-4)(k-5).
\end{equation}
It follows from Equation~(\ref{Eq-0}) that
\begin{equation}\label{Eq-B}
(q-2)(q-3)\left| PSL (2,q)_{B} \right| n \leq k(k-1)(k-2)(k-3).
\end{equation}
If we assume that $k \geq 21$, then obviously
\[k(k-1)(k-2)(k-3) < 2[(k-4)(k-5)]^2,\]
and hence
\[(q-2)(q-3)\left| PSL (2,q)_{B} \right| n< 2(q-4)^2\]
in view of Inequality~(\ref{Eq-A}). Clearly, this is only possible
when $\left| PSL (2,q)_{B} \right| \cdot n=1$. In particular, $q$
has to be even. But then the right hand side of
Equation~(\ref{Eq-0}) is always divisible by $16$ but never the left
hand side, a contradiction. If $k<21$, then the few remaining
possibilities for $k$ can easily be ruled out by hand using
Equation~(\ref{Eq-0}), Inequality~(\ref{Eq-A}), and
Corollary~\ref{divCond}.

\smallskip

\emph{Now, let us assume that $N<G \leq \aut(N)$.} We
recall that $q=p^e \geq 7$, and will distinguish in the following
the cases $p>3$, $p=2$, and $p=3$.

\smallskip

\emph{First, let $p>3$.} We define $G^*=G \cap (PSL(2,q)
\rtimes \text{\footnotesize{$\langle$}} \tau_\alpha
\text{\footnotesize{$\rangle$}})$ with $\tau_\alpha \in$
Sym$(GF(p^e) \cup \{\infty\}) \cong S_v$ of order $e$ induced by the
Frobenius automorphism $\alpha : GF(p^e) \longrightarrow GF(p^e),\,
x \mapsto x^p$. Then, by Dedekind's law, we can write
\begin{equation}\label{Eq_G^*0}
G^*= PSL(2,q) \rtimes (G^* \cap \text{\footnotesize{$\langle$}}
\tau_\alpha \text{\footnotesize{$\rangle$}}).
\end{equation}
Defining $P \Si L(2,q)= PSL(2,q) \rtimes
\text{\footnotesize{$\langle$}} \tau_\alpha
\text{\footnotesize{$\rangle$}}$, it can easily be calculated that
$P \Si L(2,q)_{0,1,\infty} = \text{\footnotesize{$\langle$}}
\tau_\alpha \text{\footnotesize{$\rangle$}}$, and
$\text{\footnotesize{$\langle$}} \tau_\alpha
\text{\footnotesize{$\rangle$}}$ has precisely $p+1$ distinct fixed
points (cf., e.g.,~\cite[Ch.\,6.4,\,Lemma\,2]{demb68}). As $p>3$, we
conclude therefore that $G^* \cap \text{\footnotesize{$\langle$}}
\tau_\alpha \text{\footnotesize{$\rangle$}} \leq G^*_B$ for some
appropriate, unique block $B \in \B$ by the definition of Steiner
\mbox{$6$-designs}. Furthermore, clearly $PSL(2,q) \cap (G^* \cap
\text{\footnotesize{$\langle$}} \tau_\alpha
\text{\footnotesize{$\rangle$}}) =1.$ Hence, we have
\begin{equation}\label{Eq_G^*}
\begin{split}
\left| B^{G^*} \right| &=  \big[G^* : G^*_{B} \big] \\
                           &= \big[PSL(2,q) \rtimes (G^* \cap \text{\footnotesize{$\langle$}}
\tau_\alpha \text{\footnotesize{$\rangle$}}) : PSL(2,q)_{B} \rtimes
(G^* \cap \text{\footnotesize{$\langle$}} \tau_\alpha
\text{\footnotesize{$\rangle$}})\big] \\
                           &= \big[PSL(2,q) : PSL(2,q)_{B} \big] \\
                           & =  \left| B^{PSL(2,q)}\right|.\\
\end{split}
\end{equation}
Thus, if we assume that \mbox{$G^* \leq \Aut(\D)$} acts already
block-transitively on $\D$, then we obtain $\left| B^{G^*}
\right|=\left| B^{PSL(2,q)}\right|=b$ in view of
Remark~\ref{equa_t=5}. Hence, $PSL(2,q)$ must also act
block-transitively on $\D$, and we may proceed as in the case when
$N=G$. Therefore, let us assume that \mbox{$G^* \leq \Aut(\D)$} does
not act block-transitively on $\D$. Then, we conclude that $\big[G :
G^* \big]=2$ and $G^*$ has exactly two orbits of equal length on the
set of blocks. Thus, by Equation~(\ref{Eq_G^*}), we obtain for the
orbit containing the block $B$ that $\left| B^{G^*} \right|=\left|
B^{PSL(2,q)}\right|=\frac{b}{2}$. As it is well-known the normalizer
of $PSL(2,q)$ in Sym$(X)$ is $P \G L (2,q)$, and hence in particular
$PSL(2,q)$ is normal in $G$. It follows therefore that we have under
$PSL(2,q)$ also precisely one further orbit of equal length on the
set of blocks. Then, proceeding similarly to the case $N=G$ for each
orbit on the set of blocks, we have (representative for the orbit
containing the block $B$) that
\begin{equation}\label{Eq-0_N<G}
\frac{(q-2)(q-3)(q-4) \left| PSL (2,q)_{B} \right| n}{2}  =
k(k-1)(k-2)(k-3)(k-4)(k-5),
\end{equation}
which gives
\begin{equation}\label{Eq-0_N<G-equiv}
(q-2)(q-3)(q-4) \left| PSL (2,q)_{B} \right| =
k(k-1)(k-2)(k-3)(k-4)(k-5),
\end{equation}
as here $n=2$. Using again
\begin{equation}
q-4 \geq (k-4)(k-5),
\end{equation}
we obtain
\begin{equation}
(q-2)(q-3) \left| PSL (2,q)_{B} \right| \leq k(k-1)(k-2)(k-3).
\end{equation}
If we assume that $k \geq 21$, then again
\[k(k-1)(k-2)(k-3) < 2[(k-4)(k-5)]^2,\]
and thus
\[(q-2)(q-3)\left| PSL (2,q)_{B} \right| < 2(q-4)^2,\]
which is only possible when $\left| PSL (2,q)_{B} \right|=1$. But,
involutions in $PSL(2,q)$ have precisely two fixed points on the
points of the projective line for \mbox{$q \equiv 1$ (mod $4$)} and
are fixed point free for \mbox{$q \equiv 3$ (mod $4$)}. Hence, each
involution always fixes a unique block by the definition of Steiner
\mbox{$6$-designs}, a contradiction. The few remaining possibilities
for $k<21$ can again easily be ruled out by hand.

\smallskip

\emph{Now, let $p=2$.} Then, clearly
$N=PSL(2,q)=PGL(2,q)$, and we have $\Aut(N)=P\Si L(2,q)$. If we
assume that $\text{\footnotesize{$\langle$}} \tau_\alpha
\text{\footnotesize{$\rangle$}} \leq P \Si L (2,q)_{B}$ for some
appropriate, unique block $B \in \B$, then, using the terminology
of~(\ref{Eq_G^*0}), we have $G^*=G=P\Si L(2,q)$ and as clearly
$PSL(2,q) \cap \text{\footnotesize{$\langle$}} \tau_\alpha
\text{\footnotesize{$\rangle$}}=1$, we can apply
Equation~(\ref{Eq_G^*}). Thus, $PSL(2,q)$ must also be
block-transitive, which has already been considered. Therefore, we
may assume that $\text{\footnotesize{$\langle$}} \tau_\alpha
\text{\footnotesize{$\rangle$}} \nleq P\Si L(2,q)_{B}$. Let $s>2$ be
a prime divisor of $e=\left| \text{\footnotesize{$\langle$}}
\tau_\alpha \text{\footnotesize{$\rangle$}} \right|$. As the normal
subgroup $H:=(P \Si L (2,q)_{0,1,\infty})^s \leq
\text{\footnotesize{$\langle$}} \tau_\alpha
\text{\footnotesize{$\rangle$}}$ of index $s$ has precisely $p^s+1$
distinct fixed points (see,
e.g.,~\cite[Ch.\,6.4,\,Lemma\,2]{demb68}), we have $G \cap H \leq
G_{B}$ for some appropriate, unique block $B \in \B$ by the
definition of Steiner \mbox{$6$-designs}. It can then be deduced
that $e=s^u$ for some $u \in \N$, since if we assume for $G= P \Si
L(2,q)$ that there exists a further prime divisor $\overline{s}>2$
of $e$ with $\overline{s} \neq s$, then $\overline{H}:=(P \Si L
(2,q)_{0,1,\infty})^{\overline{s}} \leq
\text{\footnotesize{$\langle$}} \tau_\alpha
\text{\footnotesize{$\rangle$}}$ and $H$ are both subgroups of $P\Si
L(2,q)_{B}$ by the block-transitivity of $P \Si L (2,q)$, and hence
$\text{\footnotesize{$\langle$}} \tau_\alpha
\text{\footnotesize{$\rangle$}} \leq P\Si L(2,q)_{B}$, a
contradiction. Furthermore, as $\text{\footnotesize{$\langle$}}
\tau_\alpha \text{\footnotesize{$\rangle$}}\nleq P \Si L (2,q)_{B}$,
we may, by applying Dedekind's law, assume that
\[G_{B} = PSL(2,q) _{B} \rtimes (G \cap H).\]
Thus, by Remark~\ref{equa_t=5}, we obtain
\[(q-2)(q-3)(q-4) \left| PSL (2,q)_{B} \right| \left| G \cap H \right| =
k(k-1)(k-2)(k-3)(k-4)(k-5) \left| G \cap
\text{\footnotesize{$\langle$}} \tau_\alpha
\text{\footnotesize{$\rangle$}}\right|.\] More precisely:
\begin{enumerate}
\item[(A)] if $G= PSL(2,q) \rtimes (G \cap H)$:
\[(q-2) (q-3)(q-4)\left| PSL (2,q)_{B} \right| = k(k-1)(k-2)(k-3)(k-4)(k-5)\]

\item[(B)] if $G = P \Si L (2,q)$:
\[(q-2)(q-3)(q-4) \left| PSL (2,q)_{B} \right| =k(k-1)(k-2)(k-3)(k-4)(k-5) s.\]
\end{enumerate}
As far as condition~(A) is concerned, we may argue exactly as in the
earlier case $N=G$. Thus, only condition~(B) remains. If $e$ is a
power of $2$, then Remark~\ref{equa_t=5} gives
\[(q-2)(q-3)(q-4) \left|G_{B} \right| =
k(k-1)(k-2)(k-3)(k-4)(k-5)a\] with $a \mid e$. In particular, $a$
must divide $\left|G_{B} \right|$, and we may proceed similarly as
in the case $N=G$, yielding a contradiction.

\smallskip

\emph{The case $p=3$} may be treated, mutatis mutandis, as
the case $p=2$.

\medskip

\emph{Case} (3): $N=M_v$, $v=11,12,22,23,24$.

\smallskip

By Corollary~\ref{Cameron_t=5}, we get $k=7$ for $v=11$ or $12$, and $k=7$ or $8$ for $v=22$, $23$ or $24$, and the very
small number of cases for $k$ can easily be eliminated by hand using
Corollary~\ref{divCond} and Remark~\ref{equa_t=5}.

\medskip

\emph{Case} (4): $N=M_{11}$, $v=12$.

\smallskip

As in Case~(3), for $v=12$, we have $k=7$ in view of Corollary~\ref{Cameron_t=5},
a contradiction since no \mbox{$6$-$(12,7,1)$} design can exist by Corollary~\ref{divCond}.

\medskip

This completes the proof of the Main Theorem.

\smallskip

\begin{rem}
\emph{The cases excluded from the Main Theorem remain elusive. One
can slightly reduce the possible open cases by some sophisticated
and lengthy work on condition~(B) and the corresponding one for
$p=3$. This includes a detailed consideration of the orbit-lengths
from the action of subgroups of $PSL(2,q)$ on the points of the
projective line (cf.~\cite{Hu2007a}). More precisely, we obtain the
equality
\[(q-2)(q-3)(q-4) 6 c =k(k-1)(k-2)(k-3)(k-4)(k-5) s,\]
where $q=p^{s^u}$, $p=2$ or $3$, $s^u$ some odd prime power, $s>6c$,
$c=1$, $2$, $4$ or $5$. By Siegel's classical theorem~\cite{sieg1929} on
integral points on algebraic curves only a finite number of
solutions are possible for fixed $s$. However, with regard to the
additional arithmetical conditions that are imposed in these cases,
it seems to be very unlikely that admissible parameter sets of
Steiner \mbox{$6$-designs} can be found.}
\end{rem}


\subsection*{Acknowledgment}\hfill

\smallskip

The author would like to thank the two anonymous referees for their careful reading and valuable comments.
The author gratefully acknowledges support by the Deutsche Forschungsgemeinschaft (DFG) via a Heisenberg grant (Hu954/4).

\bibliographystyle{amsplain}
\bibliography{XbibCamPraeg}

\providecommand{\bysame}{\leavevmode\hbox to3em{\hrulefill}\thinspace}
\providecommand{\MR}{\relax\ifhmode\unskip\space\fi MR }
\providecommand{\MRhref}[2]{%
  \href{http://www.ams.org/mathscinet-getitem?mr=#1}{#2}
}
\providecommand{\href}[2]{#2}
\begin{thebibliography}{10}

\bibitem{BJL1999}
Th. Beth, D.~Jungnickel, and H.~Lenz, \emph{Design {Theory}}, Vol. {I} and
  {II}, Encyclopedia of Math. and Its Applications {\bf 69/78}, Cambridge Univ.
  Press, Cambridge, 1999.

\bibitem{Block1965}
R.~E. Block, \emph{Transitive groups of collineations on certain designs},
  Pacific J. Math. \textbf{15} (1965), 13--18.

\bibitem{Buek1968}
F.~Buekenhout, \emph{Remarques sur l'homog\'{e}n\'{e}it\'{e} des espaces
  lin\'{e}aires et des syst\`{e}mes de blocs}, Math. Z. \textbf{104} (1968),
  144--146.

\bibitem{Buek1990}
F.~Buekenhout, A.~Delandtsheer, J.~Doyen, P.~B. Kleidman, M.~W. Liebeck, and
  J.~Saxl, \emph{Linear spaces with flag-transitive automorphism groups}, Geom.
  Dedicata \textbf{36} (1990), 89--94.

\bibitem{Cam1976}
P.~J. Cameron, \emph{Parallelisms of {Complete Designs}}, London Math. Soc.
  Lecture Note Series {\bf 23}, Cambridge Univ. Press, Cambridge, 1976.

\bibitem{Cam1981}
\bysame, \emph{Finite permutation groups and finite simple groups}, Bull.
  London Math. Soc. \textbf{13} (1981), 1--22.

\bibitem{CamPrae1993}
P.~J. Cameron and C.~E. Praeger, \emph{Block-transitive $t$-designs, {II}:
  large $t$}, in: Finite Geometry and Combinatorics (Deinze 1992), ed. by F. De
  Clerck et al., London Math. Soc. Lecture Note Series {\bf 191}, {Cambridge
  Univ. Press}, {Cambridge}, 1993, 103--119.

\bibitem{crc06}
C.~J. Colbourn and J.~H. Dinitz (eds.), \emph{Handbook of {Combinatorial}
  {Designs}}, 2nd ed., CRC Press, Boca Raton, 2006.

\bibitem{Atlas1985}
J.~H. Conway, R.~T. Curtis, S.~P. Norton, R.~A. Parker, and R.~A. Wilson,
  \emph{Atlas of {Finite} {Groups}}, Clarendon Press, Oxford, 1985.

\bibitem{Del2001}
A.~Delandtsheer, \emph{Finite flag-transitive linear spaces with alternating
  socle}, in: Algebraic Combinatorics and Applications, Proc. Euroconf.
  (G\"{o}{\ss}weinstein 1999), ed. by A. Betten et al., {Springer}, {Berlin},
  2001, 79--88.

\bibitem{demb68}
P.~Dembowski, \emph{Finite {Geometries}}, Springer, Berlin, Heidelberg, New
  York, 1968; Reprint 1997.

\bibitem{Gor1982}
D.~Gorenstein, \emph{{Finite Simple Groups. An Introduction} to {Their
  Classification}}, Plenum Publishing Corp., New York, London, 1982.

\bibitem{hall86}
M.~Hall{, Jr.}, \emph{Combinatorial {Theory}}, 2nd ed., J. Wiley, New York,
  1986.

\bibitem{Hu2001}
M.~Huber, \emph{Classification of flag-transitive {Steiner} quadruple systems},
  J. Combin. Theory, Series A \textbf{94} (2001), 180--190.

\bibitem{Hu2005}
\bysame, \emph{The classification of flag-transitive {Steiner} $3$-designs},
  Adv. Geom. \textbf{5} (2005), 195--221.

\bibitem{Hu_Habil2005}
\bysame, \emph{On {Highly} {Symmetric} {Combinatorial} {Designs}},
  Habilitationsschrift, Univ. T\"{u}bingen (2005), Shaker Verlag, Aachen, 2006.

\bibitem{Hu2007}
\bysame, \emph{A census of highly symmetric combinatorial designs}, J. Algebr.
  Comb. \textbf{26} (2007), 453--476.

\bibitem{Hu2007a}
\bysame, \emph{The classification of flag-transitive {Steiner} $4$-designs}, J.
  Algebr. Comb. \textbf{26} (2007), 183--207.

\bibitem{Hu_mmics2008}
\bysame, \emph{Steiner $t$-designs for large $t$}, In: Math. Methods in Comp.
  Science (MMICS) 2008, ed. by J.~Calmet et al., Lecture Notes in Comp. Science
  {\bf 5393} (Beth Festschrift), Springer, {Berlin, Heidelberg, New York},
  2008, 18--26.

\bibitem{Hu2008}
\bysame, \emph{Flag-transitive {S}teiner {D}esigns}, {Birkh\"{a}user}, Basel,
  Berlin, Boston, 2009.

\bibitem{hupi85}
D.~R. Hughes and F.~C. Piper, \emph{Design {Theory}}, Cambridge Univ. Press,
  Cambridge, 1985.

\bibitem{HupI1967}
B.~Huppert, \emph{{Endliche} {Gruppen} {I}}, Springer, Berlin, Heidelberg, New
  York, 1967.

\bibitem{Kant1972}
W.~M. Kantor, \emph{$k$-homogeneous groups}, Math. Z. \textbf{124} (1972),
  261--265.

\bibitem{Kant1985}
\bysame, \emph{Homogeneous designs and geometric lattices}, J. Combin. Theory,
  Series A \textbf{38} (1985), 66--74.

\bibitem{Kor1992}
E.~I. Khukhro and V.~D. Mazurov (eds.), \emph{Unsolved {Problems} in {Group}
  {Theory}. {The} {Kourovka} {Notebook}}, 12th ed., Russian Academy of Science,
  Novosibirsk, 1992.

\bibitem{Kleid1990}
P.~B. Kleidman, \emph{The finite flag-transitive linear spaces with an
  exceptional automorphism group}, in: Finite Geometries and Combinatorial
  Designs (Lincoln, NE, 1987), ed. by E. S. Kramer and S. S. Magliveras,
  Contemp. Math. {\bf 111}, {Amer. Math. Soc.}, {Providence, RI}, 1990,
  117--136.

\bibitem{KlLi1990}
P.~B. Kleidman and M.~W. Liebeck, \emph{The {Subgroup Structure} of the {Finite
  Classical Groups}}, London Math. Soc. Lecture Note Series {\bf 129},
  Cambridge Univ. Press, Cambridge, 1990.

\bibitem{Lieb1987}
M.~W. Liebeck, \emph{The affine permutation groups of rank three}, Proc. London
  Math. Soc. (3) \textbf{54} (1987), 477--516.

\bibitem{Lieb1998}
\bysame, \emph{The classification of finite linear spaces with flag-transitive
  automorphism groups of affine type}, J. Combin. Theory, Series A \textbf{84}
  (1998), 196--235.

\bibitem{LivWag1965}
D.~Livingstone and A.~Wagner, \emph{Transitivity of finite permutation groups
  on unordered sets}, Math. Z. \textbf{90} (1965), 393--403.

\bibitem{Luene1965}
H.~L\"{u}neburg, \emph{Fahnenhomogene {Quadrupelsysteme}}, Math. Z. \textbf{89}
  (1965), 82--90.

\bibitem{Ray-ChWil1975}
D.~K. Ray-Chaudhuri and R.~M. Wilson, \emph{On $t$-designs}, Osaka J. Math.
  \textbf{12} (1975), 737--744.

\bibitem{Saxl2002}
J.~Saxl, \emph{On finite linear spaces with almost simple flag-transitive
  automorphism groups}, J. Combin. Theory, Series A \textbf{100} (2002),
  322--348.

\bibitem{sieg1929}
C.~L. Siegel, \emph{{\"{U}}ber einige {A}nwendungen diophantischer
  {A}pproximationen}, Abh. Preuss. Akad. Wiss., Phys. Math. Kl. (1929), 41--69.

\bibitem{stin04}
D.~R. Stinson, \emph{Combinatorial {Designs:} {Constructions} and {Analysis}},
  Springer, Berlin, Heidelberg, New York, 2004.

\bibitem{Tits1964}
J.~Tits, \emph{Sur les systèmes de {Steiner} associés aux trois ``grands''
  groupes de {Mathieu}}, Rendic. Math. \textbf{23} (1964), 166--184.

\end{thebibliography}
\end{document}